\documentclass[12pt]{article}

\usepackage{amsmath}
\usepackage{amsthm}
\usepackage{amsfonts}
\usepackage{amssymb}
\usepackage{amsmath,amscd}

\newcommand{\C}{\mathbb{C}}

\newcommand{\R}{\mathbb{R}}

\newcommand{\N}{\mathbb{N}}
\newcommand{\PP}{\mathbb{P}}
\newcommand{\RR}{\mathbb{R^{+}}}

\newcommand{\codim}{\operatorname{codim}}

\newcommand{\Pic}{\operatorname{Pic}}
\newcommand{\NE}{\operatorname{NE}}
\newcommand{\NNE}{\overline{\operatorname{NE}}}

\newcommand{\Exc}{\operatorname{Exc}}

\newcommand{\OO}{{\cal O}}

\newcommand{\wtilde}{\widetilde}

\newtheorem{theo}{Theorem}
\newtheorem{prop}{Proposition}
\newtheorem{lem}{Lemma}

\title
{On weak Fano manifolds with small contractions 
obtained by blow-ups of a product of projective spaces}
\author{Toru Tsukioka}

\begin{document}

\maketitle

\begin{abstract}
We consider weak Fano manifolds with small contractions 
obtained by blowing up successively curves and subvarieties of codimension 2 in products of projective spaces. 
We give a classification result for a special case. 
In the process of proof, we describe explicitly the structure of nef cones and compute the self intersection numbers of anti-canonical divisors for such weak Fano manifolds.

\


\end{abstract}

\section{Introduction}
A smooth projective variety is called {\it Fano manifold} 
if its anti-canonical divisor is ample. 
The classification is known up to dimension 3. However, in dimension greater than or equal to 4, 
there exist only partial classification results
(see \cite{Casagrande4fold} for a recent progress). 

It is essential to investigate Fano manifolds in terms of 
the theory of extremal contractions (see \cite{MM},\cite{Wis}). 
Recall that a {\it small contraction} is a birational morphism 
whose exceptional locus has codimension greater than or 
equal to 2, and it does not appear as extremal contraction 
for smooth 3-folds. 
Hence, in dimension greater than or equal to $4$, it is interesting 
to give examples of Fano manifolds having small contractions. 

We can construct a smooth projective variety with 
a small contraction by means of successive blow-ups
(see \cite{Kawamata}): 
Let $Y$ be a smooth projective variety of dimension 
greater than or equal to 4. Let $C$ be a smooth curve on $Y$ 
and $S$ a smooth subvariety of $Y$ 
with $\codim_{Y}S=2$. Assume that $C$ and $S$ intersect 
transversally at points. Let $\pi:X\to Y$ be the blow-up 
along $C$ and let $S'$ be the strict transform of $S$ by $\pi$.
Let $\beta:\wtilde{X}\to X$ be the blow-up along $S'$. Then 
$\wtilde{X}$ has a small contraction  
(see Section \ref{small} for details). We consider the following: 

\

{\it Problem.} Classify the triples $(Y,C,S)$ 
such that $\wtilde{X}$ is a Fano manifold. 

\

The purpose of this paper is to give a classification 
result in a special case for the problem expanded to the case where 
$\wtilde{X}$ is a {\it weak Fano manifold}, i.e. 
a smooth projective variety with nef and big anti-canonical divisor. 

Throughout the paper, we work over the field of complex numbers.

\begin{theo}\label{main} Let $Y=\PP^{n-1}\times \PP^{1}$ with $n\geq 3$. 
Let $C$ be a fiber of the projection 
$Y\to \PP^{n-1}$ and let $S$ be a complete intersection of 
two divisors of bidegrees $(a,b)$ and $(1,1)$. Assume that $S$ is smooth and irreducible. Assume also that $S$ and $C$ intersect transversally at one point.  Let $\pi: X\to Y$ be the blow-up along $C$ and let 
$\beta: \wtilde{X}\to X$ be the blow-up along the strict transform of 
$S$ by $\pi$. Then 
$\wtilde{X}$ is a weak Fano manifold 
if and only if 
$n\geq 3$ and  
$$
(a,b)=(0,1),(1,0),(1,1),(1,2),(2,0),(2,1),(2,2),(3,0),(3,1) \mbox{\ or\ } (3,2).
$$
Moreover,
$\wtilde{X}$ is a Fano manifold if and only if 
$n\geq 4$ and 
$$
(a,b)=(0,1),(1,0),(1,1),(2,0) \mbox{\ or\ } (2,1).
$$
\end{theo}

\

{\it Remarks}:
(1) The case $Y=\PP^{n}$ seems more complicated 
 (see Section \ref{examplesection}). 

(2) The assumption on $C$ is not so restrictive. 
Indeed, if $C$ is not a fiber of the projection 
$p:Y=\PP^{n-1}\times \PP^{1}\to \PP^{n-1}$, 
there exists a fiber $\Gamma$ of $p$ such that 
$C\cap \Gamma\neq \emptyset$. Then we have 
$-K_{\wtilde{X}}\cdot \wtilde{\Gamma}=4-n$, 
$\wtilde{\Gamma}$ being the strict transform of $\Gamma$ by 
$\pi\circ \beta$. Hence, $-K_{\wtilde{X}}$ is not nef for $n\geq 5$.

(3) Let $q:Y=\PP^{n-1}\times \PP^{1}\to \PP^{1}$ be 
the projection. Put $y_{0}:=C\cap S$. 
Since we assume $S$ to be irreducible, $a=0$ implies 
$b=1$ and $S$ is a hyperplane in the fiber 
$q^{-1}(q(y_{0}))\simeq \PP^{n-1}$. If $a\geq 1$, 
then $q|_{S}:S\to \PP^{1}$ is surjective. 
The assumption that $S$ is contained in a divisor 
of bidegree (1,1) is natural 
(at least for the case where $\wtilde{X}$ is a Fano manifold):  Consider the open set 
$$
T:=\{t\in\PP^{1} \mid t\neq q(y_{0}) \mbox{ and } 
S\cap q^{-1}(t) \mbox{ is smooth }\}. 
$$
If $\wtilde{X}$ is a Fano manifold, 
so is $\wtilde{X}_{t}:=(q\circ \pi\circ \beta)^{-1}(t)$ for $t\in T$. 
Note that $(\pi\circ \beta)|_{\wtilde{X}_{t}}: 
\wtilde{X}_{t}\to q^{-1}(t)\simeq \PP^{n-1}$ 
is the blow-up whose center consists of the point 
$C\cap q^{-1}(t)$ and the subvariety 
$S_{t}:=S\cap q^{-1}(t)$. 
According to \cite{BCW}, there exist 
a hypersurface $U_{t}\subset q^{-1}(t)\simeq \PP^{n-1}$ of degree $a$ ($1\leq a\leq n$) 
and a hyperplane $V_{t}\subset q^{-1}(t)\simeq \PP^{n-1}$ such that 
$S_{t}$ is complete intersection of $U_{t}$ and $V_{t}$. 
Let $V$ be the closure of the union $\bigcup_{t\in T}V_{t}$. 
Then, $V$ contains $S$ and $V$ has bidegree $(1,c)$ 
for some $c\geq 0$, and our theorem covers the case $c=1$. 

\

The present paper is organized as follows:  
In Section \ref{small}, we explain how to obtain 
a small contraction by means of blow-ups. 
We also fix notations which will be used constantly 
throughout the paper. 
Section \ref{nefsection} is devoted to determine 
the structure of the nef cones of $\wtilde{X}$ 
for $(a,b)=(1,0)$ and for any $(a,b)$ such that 
$a\geq 1$ and $b\geq 0$. 
Recently, the explicite descriptions of nef cones are 
of great importance in the study of Mori dream spaces 
(see \cite{Ottem}). Hence, this section is of 
independent interest. 
In Section \ref{selfintersectionsection}, we compute 
$(-K_{\wtilde{X}})^{n}$ and express it as a 
rational function depending on $(n,a,b)$. 
We will give a sufficient condition for  
$(-K_{\wtilde{X}})^{n}$ to be strictly positive. Since the self intersection number of 
 the anti-canonical divisor is an important invariant for 
 (weak) Fano manifolds, we believe that this section is 
 also of independent interest. 
 In Section \ref{proof section}, we prove 
 Theorem \ref{main} using Propositions shown in 
 Sections \ref{nefsection} and \ref{selfintersectionsection}. 
 Section \ref{examplesection} is a supplement in which 
 we give several examples for the case $Y\neq \PP^{n-1}\times \PP^{1}$. 

\

{\it Notation}.  
Let $(x_{0}:x_{1}:\cdots : x_{n-1})$ and $(s:t)$ are   
homogeneous coordinates of $\PP^{n-1}$ and $\PP^{1}$ respectively.  
Recall that a divisor $D$ on the product $\PP^{n-1}\times \PP^{1}$ 
is said to have bidegree $(a,b)$ if 
$D$ is defined by a polynomial 
$$\sum 
c_{i_{0},i_{1}\cdots, i_{n-1}, j,k} \, 
x_{0}^{i_{0}} x_{1}^{i_{1}}\cdots x_{n-1}^{i_{n-1}}
s^{j}t^{k} \ \ 
(c_{i_{0},\cdots, i_{n-1}, j,k}\in \C)
 $$
 such that $i_{0}+\cdots +i_{n-1}=a, \ j+k=b$. 
It is equivalent to say that $D$ is a member of the linear system 
$|\OO_{\PP^{n-1}\times \PP^{1}}(a,b)|$. 

For a projective variety $X$,  we denote by 
$N^{1}(X)$ (resp. $N_{1}(X)$) the set of the numerical classes of 
divisors (resp. 1-cycles) with real coefficients. 
It is known that this is a finite dimensional vector space (see \cite{Kleiman}), and its dimension denoted by $\rho(X)$ is called the {\it Picard number} of the variety $X$. The numerical equivalence class of a divisor $D$ (resp. a 1-cycle $C$) is denoted by $[D]$ (resp. $[C]$). We see that $N^{1}(X)$ and $N_{1}(X)$ are dual to each other via the bilinear form 
$N^{1}(X)\times N_{1}(X)\to \R$ defined by the intersection number : $([D],[C])\mapsto D\cdot C$.  

The {\it nef cone} $Nef(X)$ and the {\it cone of curves}  
$\NE(X)$ are defined by 
\begin{eqnarray*}
Nef(X)&:=&\{ [D]\in N^{1}(X) \mid D \mbox{ is a nef divisor} \},  \\ 
\NE(X)&:=&\{ \sum a_{i}[C_{i}]\in N_{1}(X) 
\mid C_{i} \mbox { is an irreducible curve on } X, \ a_{i}\geq 0\}. 
\end{eqnarray*}
The closure of $\NE (X)$ in $N_{1}(X)$ is denoted by $\NNE(X)$. 
The important fact is that the two cones $Nef(X)$ and $\NNE(X)$ are dual to each other (see \cite{Lazarsfeldbook} Proposition 1.4.28). 

Let $\Gamma$ be a 1-cycle on a projective variety $Y$ and let $V$ be 
a subvariety of $Y$. For a divisor $D$ on $V$, we denote by $(D\cdot \Gamma)_{V}$ the intersection number 
taken in $V$. 
Given a birational morphism $\alpha: X\to Y$, the strict transform of 
a subvariety $M\subset Y$ will be denoted by $\alpha^{-1}_{*}M$. 

\section{Construction of a small contraction}\label{small} 

We follow Example (2.6) in \cite{Kawamata}.
Let $Y$ be a smooth projective variety of dimension $n\geq 3$. 
Let $C\subset Y$ be a smooth curve 
and let $S\subset Y$ be a smooth subvariety of codimension 2. 
Assume that $C$ and $S$ intersect transversally at one point. 
Put $y_{0}:=S\cap C$.  
Let $\pi:X\to Y$ be the blow-up along $C$ 
with the exceptional divisor $E$. Note that $\pi |_{E}:E\to C$ 
is a $\PP^{n-2}$-bundle. Put $E_{0}:=\pi^{-1}(y_{0})$. 
Let $\beta:\wtilde{X}\to X$ be the blow-up along $S':=\pi^{-1}_{*}S$ 
with the exceptional divisor $F$. 
Let $f$ be a fiber of the $\PP^{1}$-bundle $\beta |_{F}:F\to S'$. 
We put $\wtilde{E}:=\beta^{-1}(E)$ and 
$\wtilde{E_{0}}:=\beta^{-1}_{*}E_{0}$. Note that $\wtilde{E_{0}}$ 
is isomorphic to $\PP^{n-2}$. 

\begin{lem}\label{contraction} 
There exists a birational morphism $\varphi: \wtilde{X}\to X_{0}$, 
$X_{0}$ being a projective variety, 
such that $\varphi(\wtilde{E_{0}})$ is a point for $n\geq 4$. 
The same holds for $n=3$, if we assume $-K_{\wtilde{X}}$ is 
nef and big. 
\end{lem}
{\it Proof.} (See also \cite{Debarrebook} Chapter 6.)  
Let $\wtilde{e_{0}}$ be a line in $\wtilde{E_{0}}\simeq \PP^{n-2}$. 
We show that $\RR[\wtilde{e_{0}}]$ is extremal in the cone 
$\NNE(\wtilde{X})$. Assume that there exist irreducible curves 
$A,B\subset \wtilde{X}$ such that 
$\wtilde{e_{0}}\equiv A+B$. Let $D$ be an ample divisor on $Y$ and 
put $\wtilde{D}:=(\pi\circ\beta)^{*}D$. 
Since $\wtilde{D}\cdot \wtilde{e_{0}}=0$, we have 
$\wtilde{D}\cdot A=\wtilde{D}\cdot B=0$, which implies that 
$A$ and $B$ are contracted by $\pi\circ\beta$. 
Assume $A\not\subset \wtilde{E}$. 
Then there exists $s\in S\setminus \{ y_{0}\}$ such that 
$A=(\pi\circ\beta)^{-1}(s)$. Since $(\pi\circ\beta)(B)$ is a point, 
$B$ is one of the following types: 
\begin{enumerate}
\item $B=(\pi\circ\beta)^{-1}(t)$ \  
$(t\in S\setminus \{ y_{0}\})$
\item $B\subset (\pi\circ\beta)^{-1}(c)$ \ 
$(c\in C\setminus \{ y_{0}\})$ 
\item $B\subset (\pi\circ\beta)^{-1}(y_{0})$
\end{enumerate}
In case 1, we have 
$\wtilde{e_{0}}\equiv A+B\equiv f+f=2f$, a contradiction. 
In case 2, we have 
$F\cdot A+F\cdot B=-1+0=-1$, 
while $F\cdot (A+B)=F\cdot \wtilde{e_{0}}=1$, a contradiction. 
In case 3, if we put $G:=F\cap \wtilde{E}$  
then we have $(\pi\circ\beta)^{-1}(y_{0})=\wtilde{E_{0}}\cup G$. 
Assume $B\subset G$. Put $G_{0}:=F\cap \wtilde{E_{0}}$. 
Note that $N_{G_{0}/G}\simeq \OO_{\PP^{n-2}}(-1)$. 
Since $F\cdot f=-1$ and $F\cdot \wtilde{e_{0}}=1$, 
we have $F|_{G}\sim -G_{0}$. Hence,  
$$
F\cdot B=F|_{G}\cdot B=(-G_{0}\cdot B)_{G}. 
$$
On the other hand, we have 
$$
F\cdot B=F\cdot \wtilde{e_{0}}-F\cdot A=1-(-1)=2. 
$$  
Hence, $(G_{0}\cdot B)_{G}=-2<0$ which implies that 
$B\subset G_{0}\subset \wtilde{E_{0}}$. 
Thus, $[B]\in \RR[\wtilde{e_{0}}]$, a contradiction 
because $\wtilde{e_{0}}\equiv A+B\equiv f +B$.  
We conclude that all the cases $1,2,3$ do not happen. 
Therefore $A\subset \wtilde{E}$. By a similar argument, 
we also have $B\subset \wtilde{E}$. 
Now, we take intersection numbers in $\wtilde{E}$: 
$$
-1=(\wtilde{E_{0}}\cdot \wtilde{e_{0}})_{\wtilde{E}}
=(\wtilde{E_{0}}\cdot A)_{\wtilde{E}}+(\wtilde{E_{0}}\cdot B)_{\wtilde{E}}
$$
which implies 
$\wtilde{E_{0}}\cdot A<0$ or $\wtilde{E_{0}}\cdot B<0$. 
Hence, $A\subset \wtilde{E_{0}}$ or $B\subset \wtilde{E_{0}}$. 
In both cases we have 
$[A]\in \RR[\wtilde{e_{0}}]$ {\it and} $[B]\in \RR[\wtilde{e_{0}}]$. 
It follows that $\RR[\wtilde{e_{0}}]$ is an extremal ray in 
$\NNE(\wtilde{X})$.

If $n\geq 4$, we have 
$K_{\wtilde{X}}\cdot \wtilde{e_{0}}=3-n<0$. 
Hence, $\RR[\wtilde{e_{0}}]$ is a $K_{\wtilde{X}}$-negative extremal ray, and we are done 
by Contraction Theorem. 

In case $n=3$, since we assume $-K_{\wtilde{X}}$ is 
nef and big, the linear system $|-mK_{\wtilde{X}}|$ defines a 
morphism for a sufficiently large $m\in \N$ 
by Base Point Free Theorem. The Stein factorization 
gives a desired contraction because we have  
$-K_{\wtilde{X}}\cdot \wtilde{e_{0}}=0$ 
(note that $\wtilde{E_{0}}=\wtilde{e_{0}}$ for $n=3$). \qed

\

From now on, we fix the following: 

\

{\bf Notation (*)}

Assume $n\geq 3$ and put $Y:=\PP^{n-1}\times \PP^{1}.$ 
Let $p:Y\to \PP^{n-1}$ and $q:Y\to \PP^{1}$ be the projections.  
Let $C$ be a fiber of $p$. Put $H:=p^{*}\OO_{\PP^{n-1}}(1)$ and 
$L:=q^{*}\OO_{\PP^{1}}(1)$. 

Consider $V\in |H+L|$ and $U\in |aH+bL|$ where 
$a$ and $b$ are non-negative integers.   
Let $S$ be the complete intersection of $U$ and $V$. 
We assume that $S$ is smooth and irreducible. 
We assume also that $C$ and $S$ intersect transversally at 
one point and put $y_{0}:=C\cap S$. 

Let $h$ be a fiber of $p$ such that $h\neq C$ and $h\cap S=\emptyset$
and let $l$ be a line in fiber of $q$ such that 
$l\cap C=\emptyset$ and $l\cap S=\emptyset$.  

Let $\pi:X\to Y$ be the blow-up along $C$. 
Put $E:=\Exc (\pi)$ and 
$E_{0}:=\pi^{-1}(y_{0})$. 
Let $e_{0}$ be a line in $E_{0}\simeq \PP^{n-2}$
and let $e$ be a line in a fiber different from $E_{0}$ 
of the $\PP^{n-2}$-bundle $\pi |_{E}:E\to C$. 
Let $H'$ and $L'$ be the pull backs of $H$ and $L$ by $\pi$. 
Let $h'$ and $l'$ be the strict transforms of $h$ and $l$ by $\pi$.  
Put $S':=\pi^{-1}_{*}S$

Let $\beta:\wtilde{X}\to X$ be the blow-up along $S'$. 
Put $F:=\Exc (\beta)$ and $\wtilde{E}:=\beta^{-1}(E)$.  
Let $\wtilde{H}$ and $\wtilde{L}$ be the pull backs by $\beta$ of 
$H'$ and $L'$. 
Let $f$ be a fiber of the $\PP^{1}$-bundle $\beta|_{F}:F\to S'$. 
Let 
$\wtilde{e_{0}}$, 
$\wtilde{e}$,  
$\wtilde{h}$ and 
$\wtilde{l}$ 
be the strict transforms by $\beta$ of 
$e_{0}$, $e$, $h'$ and $l'$.  
Put $V':=\pi^{-1}_{*}V$ and $\wtilde{V}:=\beta^{-1}_{*}V'$. 

\section{Structure of nef cones}\label{nefsection}

The following is useful to determine the structue of simplicial cones:   
\begin{lem}\label{cone} 
Let $(D,C)\mapsto D\cdot C$ be a bilinear form of $\R^{m}\times (\R^{m})^{*}$. Let $V$ be a cone in $\R^{m}$  and 
let $V^{*}$ be its dual cone. Assume that 
there exist  
$D_{1}, D_{2}, \cdots ,D_{m}\in V$ and 
$C_{1}, C_{2}, \cdots ,C_{m}\in V^{*}$ such that 
$D_{i}\cdot C_{j}=\delta_{ij}$ (Kronecker delta). 
Then, we have 
\begin{eqnarray*}
V=\RR D_{1}+\RR D_{2}+\cdots +\RR D_{m}, \\ 
V^{*}=\RR C_{1}+\RR C_{2}+\cdots +\RR C_{m}.  
\end{eqnarray*}
\end{lem}
{\it Proof}. Since $D_{1},\cdots, D_{m}$ are linearly independent, 
for any $D\in V$ there exist real numbers $a_{1},\cdots , a_{m}$
 such that $D=a_{1}D_{i}+\cdots +a_{m}D_{m}$. We have 
 $$
 a_{i}=(a_{1}D_{1}+ \cdots + a_{m}D_{m})\cdot C_{i}=D\cdot C_{i}\geq 0 
 $$
 for $i=1,\cdots , m$. Hence,  
 $D\in \RR D_{1}+\RR D_{2}+\cdots +\RR D_{m}. $
 The structure of $V^{*}$ is similarly determined. 
 \qed

\begin{lem}\label{nef}
Let $X$ be a smooth projective variety, 
$V$ a prime divisor on $X$
 and 
$D$ a divisor on $X$. If the divisors $D-V$ and $D|_{V}$ are nef,  
then $D$ is nef. 
\end{lem}
{\it Proof.} Let $\Gamma$ be a curve on $X$. 
If $\Gamma\not\subset V$, we have 
$
D\cdot \Gamma=(D-V)\cdot \Gamma +V\cdot \Gamma\geq 0. 
$
If $\Gamma\subset V$, then we have 
$D\cdot \Gamma=D|_{V}\cdot \Gamma\geq 0$. \qed 

\

Now, we return to our situation 
(Notation (*) in Section \ref{small}). 
\begin{prop}\label{nefcone} We have 
$$
Nef(\wtilde{X})
=\RR[\wtilde{H}]
+\RR[\wtilde{L}]
+\RR[\wtilde{H}-\wtilde{E}]
+\RR[D(a,b)], 
$$
where 
$$
D(a,b):=
\begin{cases}
\wtilde{H}+\wtilde{L}-\wtilde{E}-F & 
\ \ \mbox{for $a=0$ and $b=1$}, \\ 
2\wtilde{H}+\wtilde{L}-\wtilde{E}-F & 
\ \ \mbox{for $a=1$ and $b=0$}, \\   
2\wtilde{H}+b\wtilde{L}-\wtilde{E}-F & 
\ \ \mbox{for $a=1$ and $b\geq 1$}, \\   
a\wtilde{H}+\wtilde{L}-\wtilde{E}-F & 
\ \ \mbox{for $a\geq 2$ and $b=0$}, \\   
a\wtilde{H}+b\wtilde{L}-\wtilde{E}-F & 
\ \ \mbox{for $a\geq 2$ and $b\geq 1$}.
\end{cases}
$$
\end{prop}
{\it Proof.}
We define 1-cycles $l(a)$ and $h(b)$ on $\wtilde{X}$ by:
$$
l(a):=
\begin{cases}
\wtilde{l}-\wtilde{e_{0}}-f & (a=0) \\ 
\wtilde{l}-\wtilde{e_{0}}-2f & (a=1) \\ 
\wtilde{l}-\wtilde{e_{0}}-af & (a\geq 2)
\end{cases}
, \ \ \ 
h(b):=
\begin{cases}
\wtilde{h}-f & (b=0) \\ 
\wtilde{h}-bf & (b\geq 1).
\end{cases}
$$

{\it Claim.} For any $a\geq 0$, we have $[l(a)]\in \NE(\wtilde{X})$.

{\it Proof.} Let $t_{0}:=q(y_{0})$ and 
$t\in \PP^{1}\setminus \{ y_{0}\}$. 
Put $y_{t}:=C\cap q^{-1}(t)$. Put also $Y_{0}:=q^{-1}(t_{0})$ 
and $Y_{t}:=q^{-1}(t)$.
We define the curve $\Gamma$ as follows: 
If $a=0$, let $\Gamma$ be a line in $Y_{t}\simeq \PP^{n-1}$. 
If $a=1$, let $\Gamma$ be a line in $Y_{t}$ such that 
$y_{t}\in \Gamma$ and $S\cap \Gamma\neq \emptyset$. 
If $a\geq 2$, let $\Gamma$ be a line in $Y_{0}\simeq \PP^{n-1}$ 
such that $y_{0}\in \Gamma$ and $\Gamma \subset V$. 
For any $a\geq 0$, $\Gamma\equiv l$ in $Y$. 
Put $\Gamma':=\pi^{-1}_{*}\Gamma$ and 
$\wtilde{\Gamma}:=\beta^{-1}_{*}\Gamma'$. 
For $a=0$ and $a=1$, we have $\Gamma' +e\equiv l'$. 
This yields $\wtilde{\Gamma}+\wtilde{e}\equiv \wtilde{l}$ for 
$a=0$ (because $\Gamma'\cap S'=\emptyset$) and 
$\wtilde{\Gamma}+\wtilde{e}+f\equiv \wtilde{l}$ for $a=1$ 
(because $\Gamma'$ and $S'$ intersect transversally at one point). 
In case $a\geq 2$, we have $\Gamma'+e_{0}\equiv l'$ 
which yields 
$$
(\wtilde{\Gamma}+(a-1)f)+(\wtilde{e_{0}}+f)\equiv \wtilde{l}
$$
because $(S'\cdot \Gamma')_{V'}=a-1$ 
and $(S'\cdot e_{0})_{V'}=1$. Thus, for any $a\geq 0$, we have 
$[l(a)]=[\wtilde{\Gamma}]\in \NE (\wtilde{X}). $ \qed

\

{\it Claim.} For any $b\geq 0$, we have $[h(b)]\in \NE(\wtilde{X})$.

{\it Proof.} We define the curve $\Delta$ as follows: 
If $b=0$, let $\Delta$ be a fiber of $p|_{U}$ different from $C$. 
Note that $U$ is isomorphic to $p(U)\times \PP^{1}$ because 
$U\sim aH$. If $b\geq 1$, let $\Delta$ be a fiber of $p$ such that 
$\Delta\subset V$ and $\Delta \not\subset S$ 
($\Delta$ is a fiber of the exceptional divisor of the blow-up 
$p|_{V}:V\to \PP^{n-1}$). 
Since $\Delta \equiv h$ for any $b\geq 0$, we have 
$$
(S\cdot \Delta)_{U}=V|_{U}\cdot \Delta=V\cdot h=1 \ \ \ \mbox{for } b=0,
$$
$$
(S\cdot \Delta)_{V}=U|_{V}\cdot \Delta=U\cdot h=b \ \ \ \mbox{for } b\geq 1.
$$
Put $\wtilde{\Delta}:=(\pi\circ \beta)^{-1}_{*}\Delta$. 
Then, if $b=0$, we have $\wtilde{\Delta}+f\equiv \wtilde{h}$ 
and if $b\geq 1$, $\wtilde{\Delta}+bf\equiv \wtilde{h}$. 
Thus, for any $b\geq 0$, we get 
$[h(b)]=[\wtilde{\Delta}]\in \NE (\wtilde{X})$. \qed

\

{\it Claim}. The divisors $\wtilde{H},\wtilde{L},\wtilde{H}-\wtilde{E}$ 
and $D(a,b)$ are all nef. 

{\it Proof}. We see that $H=p^{*}\OO_{\PP^{n-1}}(1)$ and 
$L=q^{*}\OO_{\PP^{1}}(1)$ are nef. Hence, so are 
$\wtilde{H}=(\pi\circ \beta)^{*}H$ and 
$\wtilde{L}=(\pi\circ \beta)^{*}L$. 
Note that  $X$ is isomorphic to 
$Bl_{z}(\PP^{n-1})\times \PP^{1}$ 
where $z$ is the point $p(C)\in \PP^{n-1}$. 
For the blow-up $\varepsilon : Bl_{z}(\PP^{n-1})\to \PP^{n-1}$ 
the divisor $\varepsilon^{*}\OO_{\PP^{n-1}}(1)-\Exc (\varepsilon)$ 
is nef. Hence, so is its pull back by the projection 
$X\to Bl_{z}(\PP^{n-1})$, which is linearly equivalent to $H'-E$. 
Therefore, $\wtilde{H}-\wtilde{E}=\beta^{*}(H'-E)$ is also nef.  

We show that $D(a,b)$ is nef 
for $(a,b)=(0,1)$ and for any $(a,b)$ such that 
$a\geq 1$ and $b\geq 0$. 

First, we consider the case $(a,b)=(0,1)$. 
Put $H_{0}:=p^{-1}(p(S))$, $H_{0}':=\pi^{-1}_{*}H_{0}$ 
and $\wtilde{H_{0}}:=\beta^{-1}_{*}H_{0}'$. 
Note that we have $S=q^{-1}(q(y_{0}))\cap H_{0}$. 
Note also that 
$\pi |_{H_{0}'}:H_{0}'\to H_{0}\simeq \PP^{n-2}\times \PP^{1}$ 
is the blow-up along $C$ and 
$\beta |_{\wtilde{H_{0}}}:\wtilde{H_{0}}\to H_{0}'$ is an isomorphism. 
Let $L_{t}$ be a fiber of $q$ such that $y_{0}\not\in L_{t}$. 
Put $\wtilde{L_{t}}:=(\pi\circ \beta)^{-1}_{*}L_{t}$. 
We see that $\wtilde{L_{t}}\cap \wtilde{H_{0}}$ and 
$F\cap \wtilde{H_{0}}$ are both fibers of the projection
$$
(q\circ\pi\circ \beta)|_{\wtilde{H_{0}}}:\wtilde{H_{0}}\to \PP^{1}. 
$$ 
Hence, we have 
$\wtilde{L}|_{\wtilde{H_{0}}}
\sim \wtilde{L_{t}}|_{\wtilde{H_{0}}}
\sim F|_{\wtilde{H_{0}}}$. Therefore, 
$$
(\wtilde{H}+\wtilde{L}-\wtilde{E}-F)|_{\wtilde{H_{0}}}
\sim (\wtilde{H}-\wtilde{E})|_{\wtilde{H_{0}}},  
$$
which is nef. 
Since $\wtilde{H_{0}}\sim \wtilde{H}-\wtilde{E}-F$, we have 
$$
(\wtilde{H}+\wtilde{L}-\wtilde{E}-F)-\wtilde{H_{0}}\sim \wtilde{L}, 
$$
which is also nef. By Lemma \ref{nef}, we conclude that 
$D(0,1)=\wtilde{H}+\wtilde{L}-\wtilde{E}-F$ is nef on $\wtilde{X}$. 

Now, we show that $D(a,b)$ is nef for $a\geq 1$ and $b\geq 0$. 
Since $F|_{\wtilde{V}}\in \Pic (\wtilde{V})$ corresponds to 
$S'\in\Pic (V')$ via the isomorphism 
$\beta |_{\wtilde{V}}:\wtilde{V}\to V'$, 
the divisor $D(a,b)|_{\wtilde{V}}$ is identified with 
the following: 
$$
\begin{cases}
(2H'+L'-E)|_{V'}-S' & (a=1, b=0) \\ 
(2H'+bL'-E)|_{V'}-S' & (a=1, b\geq 1) \\  
(aH'+L'-E)|_{V'}-S' & (a\geq 2, b=0) \\ 
(aH'+bL'-E)|_{V'}-S' & (a\geq 2, b\geq 1).  
\end{cases}
$$

Note that $\pi |_{V'}:V'\to V$ is the blow-up at the 
point $y_{0}=S\cap C$ and the exceptional divisor is 
$E\cap V'$. Hence, we have 
\begin{eqnarray*}
S'
&\sim& (\pi |_{V'})^{*}S-E|_{V'} 
\sim (\pi |_{V'})^{*}(U|_{V})-E|_{V'} 
\sim (aH'+bL')|_{V'}-E|_{V'} \\ 
&=&
\begin{cases}
(H'-E)|_{V'} & (a=1, b=0) \\
(H'+bL'-E)|_{V'} & (a=1, b\geq 1) \\ 
(aH'-E)|_{V'} & (a\geq 2, b=0) \\ 
(aH'+bL'-E)|_{V'} & (a\geq 2, b\geq 1). 
\end{cases}
\end{eqnarray*}
Therefore, $D(a,b)|_{\wtilde{V}}$ corresponds to: 
$$
\begin{cases}
(H'+L')|_{V'} & (a=1, b=0) \\
H'|_{V'} & (a=1, b\geq 1) \\ 
L'|_{V'} & (a\geq 2, b=0) \\ 
0 & (a\geq 2, b\geq 1), 
\end{cases}
$$
which is nef in any case. 

On the other hand, 
since $\wtilde{V}\sim \wtilde{H}+\wtilde{L}-F$, we have 
$$
D(a,b)-\wtilde{V}
\sim 
\begin{cases}
\wtilde{H}-\wtilde{E} & (a=1, b=0) \\
\wtilde{H}+(b-1)\wtilde{L}-\wtilde{E} & (a=1, b\geq 1) \\ 
(a-2)\wtilde{H}+(\wtilde{H}-\wtilde{E}) & (a\geq 2, b=0) \\ 
(a-2)\wtilde{H}+(b-1)\wtilde{L}+(\wtilde{H}-\wtilde{E}) & (a\geq 2, b\geq 1).  
\end{cases}
$$ 
Recall that $\wtilde{H},\wtilde{L}$ and $\wtilde{H}-\wtilde{E}$ 
are nef. Hence, so is $D(a,b)-\wtilde{V}$. By Lemma \ref{nef}, 
we conclude that $D(a,b)$ is nef. \qed

\

We have the following table of intersection numbers. 
$$
\begin{array}{c|cccc}
&\wtilde{H}&\wtilde{L}&\wtilde{E}&F \\ 
\hline
\wtilde{l} &1&0&0&0 \\
\wtilde{h}  &0&1&0&0 \\ 
\wtilde{e_{0}}  &0&0&-1&1 \\ 
f         &0&0&0&-1
\end{array}
$$
By definition of $l(a)$, $h(b)$ and $D(a,b)$, 
for $(a,b)=(0,1)$ and for any $(a,b)$ such that 
$a\geq 1$ and $b\geq 0$, we have: 
$$
\begin{array}{c|cccc}
&\wtilde{H}&\wtilde{L}&\wtilde{H}-\wtilde{E}&D(a,b) \\ 
\hline
l(a) &1&0&0&0 \\
h(b)  &0&1&0&0 \\ 
\wtilde{e_{0}}  &0&0&1&0 \\ 
f         &0&0&0&1
\end{array}
$$
Now, the proposition follows from Lemma \ref{cone} 
because we have $\rho(\wtilde{X})=4$. \qed

\

{\it Remark.} In the proof, we have also shown that 
$$
\NNE(\wtilde{X})=\RR[l(a)]+\RR[h(b)]+\RR[\wtilde{e_{0}}]+\RR[f].
$$

\section{Self intersection numbers of anti-canonical divisors}
\label{selfintersectionsection}

The purpose of this section is to prove the following: 

\begin{prop}\label{self intersection} If $a=1$, we have 
$$
(-K_{\wtilde{X}})^{n}
=\frac{(7-b)n}{2}(n-1)^{n-1}-2(n-1)(n-2)^{n-1}+(n-3)^{n}. 
$$
If $a\neq 1$, we have 
\begin{eqnarray*}
(-K_{\wtilde{X}})^{n}
=
(n-a)^{n-1}\frac{(-3a+2+ab)n+a^{2}-ab}{(a-1)^{2}} 
+(n-1)^{n-1}\frac{(a^{2}-b)n-a+b}{(a-1)^{2}} \\ -2(n-1)(n-2)^{n-1}+(n-3)^{n}. 
\end{eqnarray*}
\end{prop}

\

We prepare some lemmas. 

\begin{lem}\label{blowup} Let $D$ be a divisor 
on a smooth projective variety $Y$ of 
dimension $n\geq 3$ and let $S$ be a smooth subvariety in $Y$ 
of codimension $r\geq 2$. Let $\mu:Z\to Y$ be the blow-up along $S$. 
Let $F$ be the exceptional divisor of $\mu$. 
Then, for $k=1,2,\cdots,n$, we have 
$$
(\mu^{*}D)^{n-k}F^{k}
=(-1)^{r-1}(D|_{S})^{n-k}s_{k-r}(N^{*}_{S/Y})
$$
where $s_{k-r}(N^{*}_{S/Y})$ denotes the Segre classes of 
the conormal bundle $N^{*}_{S/Y}$.  
\end{lem} 
{\it Proof}. We follow the notation in \cite{Fulton} Chapter 3 and Appendix B,  
i.e. for a vector space $V$, the projectivization $\PP(V)$ denotes 
the set of lines in $V$.  
Consider the $\PP^{r-1}$-bundle 
$\mu |_{F}:F=\PP(N_{S/Y})\to S$. 
Let $\OO(1)$ be the dual bundle of the tautological line bundle 
$\OO(-1)$ associated to $N_{S/Y}$.  
By a definition of Segre classes, we have 
\begin{eqnarray*}
((\mu |_{F})^{*}(D|_{S}))^{n-k}\OO(1)^{k-1}
&=& ((\mu |_{F})^{*}(D|_{S}))^{n-k}\OO(1)^{(r-1)+(k-r)} \\ 
&=& (D|_{S})^{n-k}s_{k-r}(N_{S/Y}).
\end{eqnarray*}
This yields  
\begin{eqnarray*}
(\mu^{*}D)^{n-k}F^{k}
&=& (\mu^{*}D|_{F})^{n-k}(F|_{F})^{k-1} \\ 
&=& ((\mu|_{F})^{*}(D|_{S}))^{n-k}\OO(-1)^{k-1} \\ 
&=& (-1)^{k-1}(D|_{S})^{n-k}s_{k-r}(N_{S/Y}) \\ 
&=& (-1)^{k-1}(D|_{S})^{n-k}(-1)^{k-r}s_{k-r}(N^{*}_{S/Y}) \\ 
&=& (-1)^{r-1}(D|_{S})^{n-k}s_{k-r}(N^{*}_{S/Y}).
\end{eqnarray*}
\qed

\

Now, we return to our situation (Notations (*) in Section \ref{small}). 
However, in what follows, we put $h:=H|_{S}$ and 
$l:=L|_{S}$. 
\begin{lem}\label{a+b} We have 
$$
h^{n-2}=a+b, \ h^{n-3}l=a, \ l^{2}\equiv 0.
$$
\end{lem} 
{\it Proof.} Note that $S=UV\equiv (aH+bL)(H+L)$. 
Since $H^{n}=0$, $H^{n-1}L=1$ and $L^{2}\equiv 0$, 
we obtain 
\begin{eqnarray*}
h^{n-2}
&=& H^{n-2}S= aH^{n}+(a+b)H^{n-1}L=a+b, \\   
h^{n-3}l
&=& H^{n-3}LS=aH^{n-1}L+(a+b)H^{n-2}L^{2}=a, \\ 
l^{2}
&=& L^{2}S\equiv 0.
\end{eqnarray*}
\qed

\

For $a\geq 1$, we put 
$P(m):=\sum_{i=0}^{m}a^{i}$ and   
$Q(m):=\sum_{i=0}^{m}(ia^{i-1}b+(m-i)a^{i})$.

\begin{lem}\label{segre} 
For $m=1,2,\cdots, n-2$, the $m$-th Segre classe is given by  
$$
s_{m}(N^{*}_{S/Y})
=P(m)h^{m}+Q(m)h^{m-1}l.
$$
\end{lem}
{\it Proof}. Put $u:=U|_{S}$ and $v:=V|_{S}$.  
Since $S=U\cap V$ is a complete intersection, 
we have 
$$
N_{S/Y}
=N_{U/Y}|_{S}\oplus N_{V/Y}|_{S}
=U|_{S}\oplus Y|_{S}
=u\oplus v. 
$$
Hence, $N^{*}_{S/Y}=(-u)\oplus (-v)$. 
By Whitney formula, we obtain 
$$
c(N^{*}_{S/Y})=(1-u)(1-v).
$$
By the equality $c\cdot s=1$ between 
the total Chern classe and the total Segre classe, 
we get 
$$
s(N^{*}_{S/Y})
=\frac{1}{1-u}\cdot \frac{1}{1-v}
=(1+u+u^{2}+\cdots)\cdot (1+v+v^{2}+\cdots ), 
$$
whose homogeneous part of degree $m$ equals 
$\sum_{i+j=m}u^{i}v^{j}$.  
Since $l^{2}\equiv 0$, we have 
$$
u^{i}v^{j}
=(ah+bl)^{i}(h+l)^{j}
=a^{i}h^{i+j}
+(ia^{i-1}b+ja^{i})h^{i+j-1}l. 
$$ 
Therefore,  
\begin{eqnarray*}
s_{m}(N^{*}_{S/Y})
&=&\sum_{i+j=m}u^{i}v^{j} \\ 
&=&\sum_{i+j=m}(a^{i}h^{i+j}
+(ia^{i-1}b+ja^{i})h^{i+j-1}l) \\ 
&=&(\sum_{i+j=m}a^{i})h^{m}
+(\sum_{i+j=m}(ia^{i-1}b+ja^{i}))h^{m-1}l.
\end{eqnarray*}
\qed

Put 
\begin{eqnarray*}
I_{n}&:=&\sum_{k=2}^{n}\binom{n}{k}(-1)^{k}P(k-2)n^{n-k}, \\ 
I'_{n}&:=&\sum_{k=2}^{n}\binom{n}{k}(-1)^{k}kP(k-2)n^{n-k}, \\ 
J_{n}&:=&\sum_{k=2}^{n}\binom{n}{k}(-1)^{k}Q(k-2)n^{n-k}. 
\end{eqnarray*}

\begin{lem}\label{sigma}
If $a=1$, we have 
\begin{eqnarray*}
I_{n}&=&n^{n}-(2n-1)(n-1)^{n-1}, \\ 
I'_{n}&=&n(n-1)^{n-1}, \\ 
J_{n}&=&\frac{b+1}{2}((5n-2)(n-1)^{n-1}-2n^{n}).
\end{eqnarray*}

If $a\geq 2$, we have 
\begin{eqnarray*}
I_{n}&=& \frac{(n-a)^{n}+(a-1)n^{n}-a(n-1)^{n}-a(n-1)^{n}}
{a(a-1)},\\ 
I'_{n}&=& \frac{n}{a-1}((n-1)^{n-1}-(n-a)^{n-1}),\\ 
J_{n}&=&\frac{(a+b-2ab)(n-a)-ab(a-1)n}{a^{2}(a-1)^{2}}(n-a)^{n-1} \\ 
&+&\frac{(a-1)n+(a+b-2)(n-1)}{(a-1)^{2}}(n-1)^{n-1}
-\frac{a+b}{a^{2}}n^{n}.
\end{eqnarray*}
\end{lem}
{\it Proof.} 
For $a=1$, we have 
\begin{eqnarray*}
P(k-2)=k-1, \ \ \ 
Q(k-2)=\frac{b+1}{2}(k^{2}-3k+2).  
\end{eqnarray*}
If $a\geq 2$, we put $\theta:=1/(a^{2}-a)$. Then, we have  
\begin{eqnarray*}
P(k-2)
&=&\theta (a^{k}-a),  \\ 
Q(k-2)
&=&\theta^{2} 
((a+b-2ab)a^{k}+b(a-1)ka^{k}
-a^{2}(a-1)k+a^{2}(a+b-2)). 
\end{eqnarray*}
The statement is verified by direct computations 
using the following equalities for $x=a$ and $x=1$:  
\begin{eqnarray*}
\sum_{k=2}^{n}\binom{n}{k}(-x)^{k}n^{n-k}
&=&(n-x)^{n}+(x-1)n^{n}, \\ 
\sum_{k=2}^{n}\binom{n}{k}k(-x)^{k}n^{n-k}
&=& xn^{n}-xn(n-x)^{n-1}, \\ 
\sum_{k=2}^{n}\binom{n}{k}k^{2}(-x)^{k}n^{n-k}
&=& x(x-1)n^{2}(n-x)^{n-2}+xn^{n}. 
\end{eqnarray*}
\qed

\

{\it Proof of Proposition \ref{self intersection}}. 
First, we consider the case $(a,b)=(0,1)$. 
Put $L_{0}:=q^{-1}(q(y_{0}))$ and $H_{0}:=p^{-1}(p(S))$. 
Note that $S$ is a hyperplane in $L_{0}\simeq \PP^{n-1}$. 
Let $L'_{0}$ and $H'_{0}$ be the strict transforms by $\pi$ of 
$L_{0}$ and $H_{0}$. Then $S'$ is the complete intersection 
of $L'_{0}$ and $H'_{0}$. Since $H'_{0}\sim H'-E$,  
$L'_{0}\sim L'$ and $L'|_{S'}\sim 0$, we have 
$N_{S'/X}\simeq \OO_{S'}(H'-E)\oplus \OO_{S'}$. 
As in the proof of Lemma \ref{segre}, this yields 
$$
s_{m}(N^{*}_{S'/X})=(H'|_{S'}-E|_{S'})^{m} 
\ \ \mbox{for} \ m=1,2,\cdots, n-2. 
$$
On the other hand, we have 
$$
-K_{X}|_{S'}
\sim (nH'+2L'-(n-2)E)|_{S'}\sim nH'|_{S'}-(n-2)E|_{S'}. 
$$
We observe that $\pi |_{S'}: S'\to S\simeq \PP^{n-2}$ 
is the blow-up at $y_{0}$. 
We have $H'|_{S'}\sim (\pi |_{S'})^{*}\OO_{\PP^{n-2}}(1)$ 
and $\Exc (\pi |_{S'})=E|_{S'}$. 
Note that 
$$
(H'|_{S'})(E|_{S'})\equiv 0, \ 
(H'|_{S'})^{n-2}=1 \ \mbox{and} \  
(E|_{S'})^{n-2}=(-1)^{n-3}.
$$ 
By Lemma \ref{blowup} for $r=2$, 
we obtain $\beta^{*}(-K_{X})^{n-1}F=0$ 
and for $k=2,\cdots, n$,  
\begin{eqnarray*}
\beta^{*}(-K_{X})^{n-k}F^{k}
&=&-(-K_{X}|_{S'})^{n-k}s_{k-2}(N^{*}_{S'/X}) \\ 
&=&-(nH'|_{S'}-(n-2)E|_{S'})^{n-k}(H'|_{S'}-E|_{S'})^{k-2} \\ 
&=&(n-2)^{n-k}-n^{n-k}. 
\end{eqnarray*}
Since $X$ is isomorphic to $\PP^{1}\times Bl_{z}(\PP^{n-1})$ 
where $z$ is a point in $\PP^{n-1}$, we have 
$$
(-K_{X})^{n}=2n(n^{n-1}-(n-2)^{n-1}).
$$ 
It follows that 
\begin{eqnarray*}
(-K_{\wtilde{X}})^{n}
&=& (\beta^{*}(-K_{X})-F)^{n} \\ 
&=& (-K_{X})^{n}
+\sum_{k=2}^{n}\binom{n}{k}(-1)^{k}\beta^{*}(-K_{X})^{n-k}F^{k}\\ 
&=& 2n(n^{n-1}-(n-2)^{n-1})
+\sum_{k=2}^{n}\binom{n}{k}(-1)^{k}((n-2)^{n-k}-n^{n-k}) \\ 
&=&2n^{n}-(n-1)^{n}-2(n-1)(n-2)^{n-1}+(n-3)^{n}. 
\end{eqnarray*}

In case $(a,b)\neq (0,1)$, 
since we cannot necessarily describe $S'\subset X$ as a complete intersection (remark that $U'\cap V'=S'\cup E_{0}$ where $U'$ and $V'$ are the 
strict transforms by $\pi$ of $U$ and $V$), it seems hard to compute $(-K_{\wtilde{X}})^{n}$ 
directly from $(-K_{X})^{n}$. We avoid this difficulty by 
considering a flip of $\wtilde{X}$: 

{\it Step 1.} Let $\mu: Z\to Y=\PP^{n-1}\times \PP^{1}$ be the blow-up along $S=U\cap V$.  
Let $F_{Z}$ be the exceptional divisor of $\mu$. 
We have   
$$
\mu^{*}(-K_{Y})^{n}
=(-K_{Y})^{n}
=(nH+2L)^{n}=2n^{n}.
$$
By Lemma \ref{blowup} for $r=2$, we have 
$\mu^{*}(-K_{Y})^{n-1}F_{Z}=0$ and 
$$
\mu^{*}(-K_{Y})^{n-k}F_{Z}^{k}
=-(-K_{Y}|_{S})^{n-k}s_{k-2}(N^{*}_{S/Y})
 \ \ \mbox{for} \ k=2,\cdots, n.
$$
Therefore, 
\begin{eqnarray*}
(-K_{Z})^{n}
&=& (\mu^{*}(-K_{Y})-F_{Z})^{n} \\ 
&=& \mu^{*}(-K_{Y})^{n}
+\sum_{k=1}^{n}\binom{n}{k}\mu^{*}(-K_{Y})^{n-k}(-F_{Z})^{k} \\ 
&=& 2n^{n}-\sum_{k=2}^{n}\binom{n}{k}
(-1)^{k}(-K_{Y}|_{S})^{n-k}s_{k-2}(N^{*}_{S/Y}). 
\end{eqnarray*}
Here, we have 
$
-K_{Y}|_{S}\sim (nH+2L)|_{S}=nh+2l. 
$
By Lemma \ref{a+b} and \ref{segre},  
$(-K_{Y}|_{S})^{n-k}s_{k-2}(N^{*}_{S/Y})$ is equal to  
$$
(3a+b)P(k-2)n^{n-k}
-2aP(k-2)kn^{n-k-1}+aQ(k-2)n^{n-k}. 
$$
Thus, 
\begin{eqnarray*}
(-K_{Z})^{n}
=2n^{n}-(3a+b)I_{n} +\frac{2a}{n}I'_{n} -aJ_{n}. 
\end{eqnarray*}

Now, substituting the result of Lemma \ref{sigma}, we conclude that  
if $a=1$, we have 
$$
(-K_{Z})^{n}
=\frac{(7-b)n}{2}(n-1)^{n-1}, 
$$
if $a\geq 2$, we have 
\begin{eqnarray*}
(-K_{Z})^{n}
=
(n-a)^{n-1}\frac{(-3a+2+ab)n+a^{2}-ab}{(a-1)^{2}} 
+(n-1)^{n-1}\frac{(a^{2}-b)n-a+b}{(a-1)^{2}}. 
\end{eqnarray*}

\

{\it Step 2.} 
Let $\alpha:\wtilde{Z}\to Z$ be the blow-up along 
the curve $C':=\mu^{-1}_{*}C$ 
with the exceptional divisor $G$. 
Note that 
$N_{C'/Z}=\OO_{\PP^{1}}^{\oplus (n-2)}\oplus \OO_{\PP^{1}}(-1)$.
Hence, we have $-K_{Z}\cdot C'=1$ 
and $\deg (N^{*}_{C'/Z})=1$. 
Using Lemma \ref{blowup} for $r=n-1$, we have the following: 
\begin{eqnarray*}
\alpha^{*}(-K_{Z})^{n-k}G^{k}
&=&0 \ \ \mbox{for} \ k=1,2,\cdots, n-2,  \\  
\alpha^{*}(-K_{Z})G^{n-1}
&=&(-1)^{n}(-K_{Z}\cdot C')=(-1)^{n}, \\ 
G^{n}&=&(-1)^{n}s_{1}(N^{*}_{C'/Z})=(-1)^{n+1}.
\end{eqnarray*} 
Therefore, 
\begin{eqnarray*}
(-K_{\wtilde{Z}})^{n}
&=&(\alpha^{*}(-K_{Z})-(n-2)G)^{n} \\ 
&=&(\alpha^{*}(-K_{Z}))^{n}
+(-1)^{n-1}n(n-2)^{n-1}\alpha^{*}(-K_{Y})G^{n-1}
+(-1)^{n}(n-2)^{n}G^{n} \\ 
&=&(-K_{Z})^{n}-n(n-2)^{n-1}-(n-2)^{n} \\ 
&=&(-K_{Z})^{n}-2(n-1)(n-2)^{n-1}. 
\end{eqnarray*}

{\it Step 3.} 
We observe that $\wtilde{X}$ and $\wtilde{Z}$ are connected by a flip. 
We have  
$$
(-K_{\wtilde{X}})^{n}=(-K_{\wtilde{Z}})^{n}+(n-3)^{n}. 
$$
To see this, put $\Gamma_{0}:=\mu^{-1}(y_{0})$ and 
$\wtilde{\Gamma_{0}}:=\alpha^{-1}_{*}\Gamma_{0}$. 
Let $\gamma :W\to \wtilde{Z}$ be the blow-up along 
the curve $\wtilde{\Gamma_{0}}$ and 
let $M$ be the exceptional divisor of $\gamma$. 
Note that $M\simeq \PP^{n-2}\times \PP^{1}$ and 
$N_{M/W}\simeq \OO_{\PP^{n-2}\times \PP^{1}}(-1,-1)$. 
The contraction map sending $M$ to $\PP^{n-2}$ is nothing but 
the blow-up $\delta :W\to \wtilde{X}$ along $\wtilde{E_{0}}\simeq \PP^{n-2}$. 

We have 
$K_{W}\sim \delta^{*}K_{\wtilde{X}}+M$ and 
$K_{W}\sim \gamma^{*}K_{\wtilde{Z}}+(n-2)M$. 
Hence, 
$$
\delta^{*}(-K_{\wtilde{X}})
\sim\gamma^{*}(-K_{\wtilde{Z}})-(n-3)M. 
$$
Note that 
$N_{\wtilde{\Gamma_{0}}/\wtilde{Z}}\simeq \OO_{\PP^{n-1}}(-1)^{\oplus (n-1)}$. 
Hence, we have $-K_{\wtilde{Z}}\cdot \wtilde{\Gamma_{0}}=3-n$ and 
$\deg(N^{*}_{\wtilde{\Gamma_{0}}/\wtilde{Z}})=n-1$. As in Step 2,  
we obtain $\gamma^{*}(-K_{\wtilde{Z}})M^{n-1}=(-1)^{n+1}(n-3)$, 
$M^{n}=(-1)^{n+1}(n-1)$ and 
$\gamma^{*}(-K_{\wtilde{Z}})^{n-k}M^{k}=0$ for $k=1,\cdots,n-2$. 

Thus, 
\begin{eqnarray*}
(-K_{\wtilde{X}})^{n}
&=&(\delta^{*}(-K_{\wtilde{X}}))^{n} \\ 
&=&(\gamma^{*}(-K_{\wtilde{Z}})-(n-3)M)^{n} \\ 
&=&(\gamma^{*}(-K_{\wtilde{Z}}))^{n}
+(-1)^{n-1}n(n-3)^{n-1}\gamma^{*}(-K_{\wtilde{Z}})M^{n-1}
+(-1)^{n}(n-3)^{n}M^{n} \\ 
&=& (-K_{\wtilde{Z}})^{n}+(n-3)^{n}.  
\end{eqnarray*}

\

By Step 2 and Step 3, we obtain
$$
(-K_{\wtilde{X}})=(-K_{Z})^{n}-2(n-1)(n-2)^{n-1}+(n-3)^{n}.
$$
Substituting the result of Step 1, we complete the proof of Proposition \ref{self intersection}. 
\qed

\

\begin{prop}\label{positive}
For $(a,b)=(0,1)$ and for $(a,b)$ such that 
$1\leq a\leq 3$ and $0\leq b\leq 3$, 
we have $(-K_{\wtilde{X}})^{n}>0$ for any $n\geq 3$. 
\end{prop}
{\it Proof.} For each case, 
we compute $(-K_{\wtilde{X}})^{n}$ 
using Proposition \ref{self intersection}. 
Note that we have $(n-1)^{n}-2(n-1)(n-2)^{n-1}>0$ for any $n\geq 3$. 

If $(a,b)=(0,1)$, then 
\begin{eqnarray*}
(-K_{\wtilde{X}})^{n}
&=& 2n^{n}-(n-1)^{n}-2(n-1)(n-2)^{n-1}+(n-3)^{n} \\ 
&\geq& 2(n-1)^{n}-(n-1)^{n}-2(n-1)(n-2)^{n-1}+(n-3)^{n} \\ 
&\geq& (n-1)^{n}-2(n-1)(n-2)^{n-1}+(n-3)^{n} \\ 
&>& 0. 
\end{eqnarray*}

Assume $b\leq 3$. 
If $a=1$,  
\begin{eqnarray*}
(-K_{\wtilde{X}})^{n}
&=&\frac{7-b}{2}n(n-1)^{n-1}-2(n-1)(n-2)^{n-1}+(n-3)^{n} \\ 
&\geq & 2n(n-1)^{n-1}-2(n-1)(n-2)^{n-1}+(n-3)^{n} \\ 
&>&(n-1)^{n}-2(n-1)(n-2)^{n-1}+(n-3)^{n} \\ 
&>&0.
\end{eqnarray*}
If $a=2$,  
\begin{eqnarray*}
(-K_{\wtilde{X}})^{n}
&=&(4n-2)(n-1)^{n-1}-6(n-1)(n-2)^{n-1} \\ 
&& -b((n-1)^{n}-2(n-1)(n-2)^{n-1})+(n-3)^{n} \\ 
&\geq & (4n-2)(n-1)^{n-1}-6(n-1)(n-2)^{n-1} \\ 
&& -3((n-1)^{n}-2(n-1)(n-2)^{n-1})+(n-3)^{n} \\ 
&=&(n-1)^{n}+(n-3)^{n} \\ 
&>& 0.
\end{eqnarray*}
If $a=3$,  
\begin{eqnarray*}
4(-K_{\wtilde{X}})^{n}
&=&-3(n+1)(n-3)^{n-1}-8(n-1)(n-2)^{n-1}+(9n-3)(n-1)^{n-1} \\ 
&&-b((n-1)^{n}-3(n-1)(n-3)^{n-1}) \\ 
&\geq & -3(n+1)(n-3)^{n-1}-8(n-1)(n-2)^{n-1}+(9n-3)(n-1)^{n-1} \\ 
&&-3((n-1)^{n}-3(n-1)(n-3)^{n-1}) \\ 
&=&(6n-12)(n-3)^{n-1}+6n(n-1)^{n-1}-8(n-1)(n-2)^{n-1} \\ 
&>& 6n(n-1)^{n-1}-8(n-1)(n-2)^{n-1} \\ 
&>& 4((n-1)^{n-1}-2(n-1)(n-2)^{n-1}) \\ 
&>&0.
\end{eqnarray*}
\qed

\

{\it Remark.} 
More precise estimations show that we have 
$(-K_{\wtilde{X}})^{n}>0$ for any $n\geq 3$ in the cases: 
$a=1$ and $b\leq 5$; $a=2$ and $b\leq 6$; $a=3$ and $b\leq 8$.
In case $a\geq 4$, 
the positivity of $(-K_{\wtilde{X}})^{n}$ is independent of 
the value $b$. For example, if $a=15$ we have 
$(-K_{\wtilde{X}})^{4}=-306b-285<0$ and 
$(-K_{\wtilde{X}})^{5}=3056b+1344>0$ for any $b\in \N$.

\section{Proof of Theorem \ref{main}}\label{proof section}

By the canonical bundle formula 
for the blow-ups $\pi$ and $\beta$, 
we have
$$
K_{\wtilde{X}}
\sim \beta^{*}K_{X}+F 
\sim \beta^{*}(\pi^{*}K_{Y})+(n-2)\wtilde{E}+F.
$$
Combining with $-K_{Y}\sim nH+2L$, we get 
$$
-K_{\wtilde{X}}\sim n\wtilde{H}+2\wtilde{L}-(n-2)\wtilde{E}-F. 
$$
First, we consider the case $a\geq 2$ and $b\geq 1$. 
We rewrite $-K_{\wtilde{X}}$ by means of the generators 
of $Nef(\wtilde{X})$: 
$$
-K_{\wtilde{X}}
=(3-a)\wtilde{H}+(2-b)\wtilde{L}
+(n-3)(\wtilde{H}-\wtilde{E})
+(a\wtilde{H}+b\wtilde{L}-\wtilde{E}-F).
$$
By Proposition \ref{nefcone}, we see that 
$-K_{\wtilde{X}}$ is nef if and only if $3-a\geq 0$, $2-b\geq 0$ and $n-3\geq 0$. Since the numerical equivalence classes of ample divisors are interior points of 
the nef cone (\cite{Lazarsfeldbook} Theorem 1.4.23), it follows that 
$-K_{\wtilde{X}}$ is ample if and only if $3-a>0$, $2-b>0$ and $n-3>0$.

In the other cases, we argue similarly in the following forms: 
\begin{itemize}
\item For $a=0$ and $b=1$, 
$$
-K_{\wtilde{X}}
\sim 
2\wtilde{H}+\wtilde{L}
+(n-3)(\wtilde{H}-\wtilde{E})+(\wtilde{H}+\wtilde{L}-\wtilde{E}-F).  
$$
\item For $a=1$ and $b=0$, 
$$
-K_{\wtilde{X}}
\sim \wtilde{H}+\wtilde{L}
+(n-3)(\wtilde{H}-\wtilde{E})+(2\wtilde{H}+\wtilde{L}-\wtilde{E}-F). 
$$
\item For $a=1$ and $b\geq 1$, 
$$
-K_{\wtilde{X}}
\sim \wtilde{H}+(2-b)\wtilde{L}
+(n-3)(\wtilde{H}-\wtilde{E})+(2\wtilde{H}+b\wtilde{L}-\wtilde{E}-F).  
$$
\item For $a\geq 2$ and $b=0$, 
$$
-K_{\wtilde{X}}
\sim 
(3-a)\wtilde{H}+\wtilde{L}
+(n-3)(\wtilde{H}-\wtilde{E})+(a\wtilde{H}+\wtilde{L}-\wtilde{E}-F).  
$$
\end{itemize}
Finally, we conclude that 
$-K_{\wtilde{X}}$ is nef if and only if 
$n\geq 3$ and  
$$
(a,b)=(0,1),(1,0),(1,1),(1,2),(2,0),(2,1),(2,2),(3,0),(3,1) \mbox{\ or\ } (3,2).
$$
Moreover, $-K_{\wtilde{X}}$ is ample if and only if 
$n\geq 4$ and 
$$
(a,b)=(0,1),(1,0),(1,1),(2,0) \mbox{\ or\ } (2,1).
$$

In general, a nef divisor $D$ is big if and only if $D^{n}>0$
(\cite{Lazarsfeldbook} Theorem 2.2.16). Thus, the proof of Theorem is completed by Proposition \ref{positive}. \qed

\section{Other examples}\label{examplesection}
Details in this section can be verified by the methods in Sections  \ref{nefsection}
 and \ref{selfintersectionsection}. We keep all the notations (*) in Section \ref{small} 
except that the divisors $H$ and $L$ are replaced by appropriate ones. 

In the following examples 1,2 and 3, we put $Y=\PP^{n}$ with $n\geq 4$. 
 
{\it Example 1}. Let $C$ be a line 
and $S$ an ($n-2$)-plane. 
Assume $C\cap S\neq \emptyset$. 
We consider $H:=\OO_{\PP^{n}}(1)$ and 
$\wtilde{H}:=(\pi\circ\beta)^{*}H$. Then, we have 
$$
Nef (\wtilde{X})
=\RR[\wtilde{H}]+\RR[\wtilde{H}-\wtilde{E}]
+\RR[2\wtilde{H}-\wtilde{E}-F],
$$
$$
-K_{\wtilde{X}}
\sim (n+1)\wtilde{H}-(n-2)\wtilde{E}-F
=2\wtilde{H}+(n-3)(\wtilde{H}-\wtilde{E})+(2\wtilde{H}-\wtilde{E}-F).
$$
Hence, $\wtilde{X}$ is a Fano manifold for any $n\geq 4$. 

{\it Example 2}. Let $C$ be a line. 
Let $S$ be the complete intersection of a hyperplane and 
a hyperquadric. Assume that $C$ and $S$ intersect transversally at 
one point. Then $\wtilde{X}$ is a Fano manifold for any $n\geq 4$. 
Indeed, the structure of nef cone and the description of the anti-canonical divisor for $\wtilde{X}$ are completely same as in Example 1. 
Even if the intersection $C\cap S$ consists of {\it two} points, 
$\wtilde{X}$ remains Fano, while the exceptional locus of the small contraction has two irreducible components. 

{\it Example 3}. Let $P\subset Y=\PP^{n}$ be a 2-plane and $C$ a smooth conic 
on $P\simeq \PP^{2}$. Let $S$ be an ($n-2$)-plane such that 
$\sharp (C\cap S)=2$. Then, we have 
$$
Nef (\wtilde{X})
=\RR[\wtilde{H}]+\RR[2\wtilde{H}-\wtilde{E}]
+\RR[3\wtilde{H}-\wtilde{E}-F],
$$
$$
-K_{\wtilde{X}}
\sim (n+1)\wtilde{H}-(n-2)\wtilde{E}-F
=(4-n)\wtilde{H}+(n-3)(2\wtilde{H}-\wtilde{E})+(3\wtilde{H}-\wtilde{E}-F).
$$
We see that $-K_{\wtilde{X}}$ is nef only for $n=4$.  Moreover, 
we have $(-K_{\wtilde{X}})^{4}=353>0$. Hence $\wtilde{X}$ is a 
weak Fano manifold for $n=4$. 

{\it Example 4}. Put $Y:=\PP^{2}\times \PP^{2}$. Let $C$ be a line 
in a fiber of a projection $Y\to \PP^{2}$ and $S$ a fiber of the other 
projection such that $C\cap S\neq\emptyset$. 
Then $\wtilde{X}$ is a Fano 4-fold. Indeed, we are able to show: 

\begin{prop}
Let $Y:=\PP^{n-2}\times \PP^{2}$ with $n\geq 3$. 
Let $C$ be a smooth plane curve of degree $d$ in a fiber of 
the projection $p:Y\to \PP^{n-2}$. Let $S$ be a 
fiber of the projection $q:Y\to \PP^{2}$ such that 
$C\cap S\neq \emptyset$. Then $\wtilde{X}$ is a weak Fano manifold 
if and only if 
$$
(n,d)=(3,1),(3,2),(3,3),(4,1) \ or \ (5,1).
$$
Moreover, $\wtilde{X}$ is a Fano manifold if and only if $(n,d)=(4,1)$. 
\end{prop}
 {\it Proof}. Put $H:=p^{*}\OO_{\PP^{n-2}}(1)$ and 
 $L:=q^{*}\OO_{\PP^{2}}(1)$. Then we have 
 $$
 Nef(\wtilde{X})=
 \RR[\wtilde{H}]+\RR[\wtilde{L}]
 +\RR[\wtilde{H}+d\wtilde{L}-\wtilde{E}]
 +\RR[\wtilde{H}+d\wtilde{L}-\wtilde{E}-F].  
 $$
 Since $K_{\wtilde{X}}
 \sim (\pi\circ \beta)^{*}K_{Y}+(n-2)\wtilde{E}+F$ and 
 $-K_{Y}\sim (n-1)H+3L$, we have 
 \begin{eqnarray*}
 -K_{\wtilde{X}}
 &\sim& (n-1)\wtilde{H}+3\wtilde{L}-(n-2)\wtilde{E}-F \\ 
 &=&\wtilde{H}+(3-d(n-2))\wtilde{L}
 +(n-3)(\wtilde{H}+d\wtilde{L}-\wtilde{E})
 +(\wtilde{H}+d\wtilde{L}-\wtilde{E}-F). 
 \end{eqnarray*}
 Hence, $-K_{\wtilde{X}}$ is nef (resp. ample) if and only if 
 $3-d(n-2)$ and $n-3$ are positive (resp. strictly positive). 
 On the other hand, we obtain
 $$
 (-K_{\wtilde{X}})^{n}
 =4n(n-1)^{n-1}+(n-2)^{n-1}(d(d-3)n-2d^{2}+2)+(n-3)^{n}, 
 $$
 which is strictly positive for 
 $(n,d)=(3,1),(3,2),(3,3),(4,1)$ and $(5,1)$. \qed 
 
 \
 
{\it Acknowledgements.}\  
The author would like to thank Kento Fujita and Kazunori Yasutake for helpful comments.

------------

{\small {\sc Department of Mathematics, 
Tokai University, 
4-1-1 Kitakaname, Hiratsuka-shi,
Kanagawa, 259-1292 Japan} }

{\it E-mail address}: \ {\tt tsukioka@tokai-u.jp}

\end{document}